 \documentclass{ajour}
\usepackage{amssymb}
\begin{document}

%







\title{On the zero modes of Pauli operators\thanks{Supported by
the European Union
under the TMR grant FMRX-CT 96-0001. } }


\authors{ A.~A.~Balinsky  and W.~D.~Evans}
\affil{School of Mathematics,
         Cardiff University,
         23 Senghennydd Road, P.O. Box~ 926,
         Cardiff~CF24~4YH,
         UK\thanks{email: BalinskyA@cardiff.ac.uk; EvansWD@cardiff.ac.uk.}}

\abstract{ Two results are proved for $\mathrm{nul} \
\mathbb{P}_A$, the dimension of the kernel of the Pauli operator
$\mathbb{P}_A = \bigl\{ \bbf{\sigma} \cdotp \bigl(\frac{1}{i}
\bbf{\nabla} + \vec{A} \bigr) \bigr\} ^2 $  in $[L^2
(\mathbb{R}^3)]^2$:
\ (i) for $|\vec{B}| \in L^{3/2}
(\mathbb{R}^3),$ where $\vec{B} = \mathrm{curl} \vec{A}$ is the
magnetic field, $\mathrm{nul} \
\mathbb{P}_{tA} = 0$ except for a finite number of values of $t$ in any
compact subset of $(0, \infty)$; \ (ii) \ $\bigl\{ \    \vec{B}:  \  \mathrm{nul} \  \mathbb{P}_{ A} = 0,
\ \ | \vec{B} | \in L^{3/2}(\mathbb{R}^3)
\ \bigr\} $ contains an open dense subset of
$[L^{3/2}(\mathbb{R}^3)]^3$. }


\authorrunninghead{A. Balinsky and W.D. Evans}
\titlerunninghead{On the zero modes of Pauli operators}


\begin{article}
\section{Introduction}

    The Pauli operator is formally defined by
\begin{equation} \label{Pauli}
\mathbb{P}_A = \biggl\{ \bbf{\sigma} \cdotp \biggl(\frac{1}{i}
\bbf{\nabla} + \vec{A} \biggr) \biggr\} ^2 \equiv \sum\limits_{j=1}^3
\biggl\{ \sigma_j  \  \biggl(\frac{1}{i} \partial_j + A_j \biggr) \biggr\} ^2
\end{equation}
where $\vec{A} = (A_1, A_2, A_3)$ is a vector potential which is such
that $\mathop{curl} \vec{A} = \vec{B}$,  the magnetic field, and
$\bbf{\sigma} = \vec{\sigma} \equiv (\sigma_1,\sigma_2,\sigma_3 )$ is the triple of Pauli
matrices
\begin{equation} \label{Paili_matr}
\sigma_1  \ = \
\left( \begin{array}{cc}
 0  & \ \  1  \\
 1  &  \ \  0
\end{array} \right) , \ \
\sigma_2  \ = \
\left( \begin{array}{cc}
 0 & \ \  -i  \\
 i  &  \ \  0
\end{array} \right) , \  \
\sigma_3  \ = \
\left( \begin{array}{cc}
 1  & \ \  0  \\
 0  &  \ \  -1
\end{array} \right).
\end{equation}
The expression (\ref{Pauli}) defines a non-negative self-adjoint
operator in $[L^2 (\mathbb{R}^3)]^2$; its precise
definition will be given in \S 2.

    Zero modes of $\mathbb{P}_A$ are the eigenvectors
corresponding to an eigenvalue at zero. The existence of  zero modes
has profound implications to the stability of matter when
$\mathbb{P}_A$, or the Dirac-Weyl operator $\bbf{\sigma} \cdotp \bigl(\frac{1}{i}
\bbf{\nabla} + \vec{A} \bigr)$, is used for the model, for the
vanishing kinetic energy of zero modes means that their potential
energy can not be controlled by their  kinetic. For an account of this
phenomenon and its consequences, we refer to \cite{FLL,LL,LY}.  Also, the
importance of zero modes for the understanding of other deep physical
problems is emphasized in \cite{AMN1}.  A significant mathematical
implication of zero modes is that there can't be an analogue of the
Cwikel-Lieb-Rosenblum inequality for the number of negative
eigenvalues of  $\mathbb{P}_A  + V$ in terms of some $L^p$ norm of the
scalar potential $V$, since any small negative perturbation $V$ would
produce negative eigenvalues, contrary to such an  inequality if $V$ is
sufficiently small.

    The first example  of a magnetic field $\vec{B}$  which yields
zero modes was the following constructed in \cite{LY}:
\[
\vec{B} (\bbf{x}) \ = \ \frac{12}{(1+ r^2)^3} \  (2 x_1 x_3 -2 x_2, \ 2 x_2 x_3 + 2
x_1,  \ 1- x_1^2 -x_2^2 + x_3^2),
\]
where $\bbf{x} = (x_1, x_2, x_3)$ and $r = |\bbf{x} |$. There are two
features of the Loss/Yau example  which are of particular relevance
to us:
 \begin{equation} \label{3}
|\vec{B} (\bbf{x}) | \in L^p (\mathbb{R}^3) \ \ \
\mbox{ for any} \ \
p > \frac{3}{4},
\end{equation}
\begin{equation} \label{4}
\frac{1}{3} (3+2l )\vec{B} (\bbf{x}), \  l \in \mathbb{N}, \ \ \
\mbox{ also yields zero modes.}
\end{equation}
We shall reserve comment on these till  later.  Other examples of zero
modes, based on the construction of \cite{LY}, are given in \cite{E},
\cite{AMN1}.
In an attempt  to explain the origin of zero modes, Erd\H{o}s and
Solovej in \cite{ES} give a more geometric viewpoint. Using the known
behaviour of the Dirac operator under conformal transformations, and
that $\mathbb{R}^3$ is conformally equivalent to a punctured sphere
$\mathbb{S}^3$, they establish their zero modes on $\mathbb{S}^3$ as
well as on $\mathbb{R}^3$ as pull-backs of zero modes on
$\mathbb{S}^2$ under the Hopf map $\mathbb{S}^3 \rightarrow
\mathbb{S}^2$.
It is also shown in \cite{ES} that arbitrary degeneracy is possible;
examples of this may also be found in \cite{AMN2}.

    In even-dimensional manifolds, the Atiyah-Singer index theorem
is a powerful tool for investigating the kernel of $\mathbb{P}_A$,
since the deficiency of $\mathbb{P}_A$ can vanish,  in which case the index
is equal to  $\mathrm{nul} \
\mathbb{P}_A$,  the nullity of $\mathbb{P}_A$  (i.e. the dimension  of the kernel $\ker
\mathbb{P}_A$).  A celebrated  example is the Aharonov-Casher Theorem in
$\mathbb{R}^2$ and its analogue due to Avron and Tomaras in
$\mathbb{S}^2$ (see \cite{CFKS}). In $\mathbb{R}^2$, this assert that for
suitable $\vec{B}$ (e.g. $\vec{B}$ bounded and of compact support),
the nullity of $\mathbb{P}_A$ is
\[
\biggl\{ \frac{1}{2 \pi}  \biggl|
\int\limits_{\mathbb{R}^2}  B (\bbf{x}) d \bbf{x} \biggr|  \biggr\}
\]
where $\{ y \}$ denotes the largest integer strictly less than $y$ and
$\{ 0 \} = 0$; note that in $\mathbb{R}^2$, the  magnetic field has
only one component, and is thus a scalar field. Thus in
$\mathbb{R}^2$, zero modes are abundant; they exist as long as the
magnetic flux $ \frac{1}{2 \pi}
\int\limits_{\mathbb{R}^2}  B (\bbf{x}) d \bbf{x} $ takes values
outside $[-1. 1]$. The situation in  $\mathbb{S}^2$ is very
different. There are now zero modes if and only if the magnetic flux is
an integer, a picture which is somewhat reminiscent of that in
(\ref{4}). In fact, this is typical of what happens on any compact
manifold (of even or odd dimension), as shown by Anghel in \cite{A}.

    Apart from the examples in \cite{LY,ES}  mentioned above very
little is known for $\mathbb{R}^3$, and indeed for non-compact
manifolds of odd dimension, since it is not easy to obtain information
from the Atiyah-Singer index theorem in this case. However, we prove that
the situation in $\mathbb{R}^3$ is like that described above for
compact manifolds, and is thus dramatically different to
$\mathbb{R}^2$.
Specifically, we prove in Theorems~4.1 and 4.2
\begin{itemize}
\item \ for $| \vec{B} | \in L^{3/2}(\mathbb{R}^3)$, \ $\mathrm{nul} \
\mathbb{P}_{t A} = 0$ except for a finite number of values of $t$  in any
compact subset of $[0, \infty)$;
\item \ $\bigl\{ \    \vec{B}:  \  \mathrm{nul} \  \mathbb{P}_{ A} = 0,
\ \mathrm{curl} \ \vec{A} = \vec{B} \ \ \mbox{and} \ \ | \vec{B} | \in L^{3/2}(\mathbb{R}^3)
\ \bigr\} $ contains an open dense subset of
$[L^{3/2}(\mathbb{R}^3)]^3$.
\end{itemize}

This explains why zero modes are so difficult to obtain. Note that the
Loss-Yau example satisfies our hypothesis.  The analogous result holds
for $\mathbb{R}^n$ with $n > 3$.

\section{Preliminaries}

    We can write (\ref{Pauli}) as
\begin{equation} \label{Paili_Zeeman}
\mathbb{P}_A =  S_A + \vec{\sigma} \cdotp \vec{B}, \ \ \  \vec{B}= \mathop{curl} \vec{A} ,
\end{equation}
where $S_A$ is the magnetic Schr\"{o}dinger  operator
\begin{equation} \label{MagneticSch}
S_A = \biggl(  \frac{1}{i} \bbf{\nabla} + \vec{A} \biggr)^2  \
\mathbb{I}_2 \ \equiv  \ \sum\limits_{j=1}^3
\biggl(\frac{1}{i} \partial_j + A_j \biggr)  ^2  \  \mathbb{I}_2,
\end{equation}
$\mathbb{I}_2$ being the $2 \times 2$ identity matrix and
$\vec{\sigma} \cdotp \vec{B}$ the Zeeman term.
Note that a gauge transformation $\vec{A} \mapsto \vec{A} + df$ does
not alter the nullity, and hence $\mathrm{nul} \ \mathbb{P}_A$ is
independent of the gauge.
We denote $[L^2
(\mathbb{R}^3)]^2$ by $\mathcal{H}$ and its standard  inner-product and
norm by $(\cdotp , \cdotp)$ and $\| \cdotp \|$ respectively:
\[ \| f \|^2 \ = \ \int\limits_{\mathbb{R}^3} \ |f(\mathbf{x})|^2
\ d  \mathbf{x},
\]
where $| \cdotp |$ is the Euclidean norm on $\mathbb{C}^2$. It will be
assumed throughout that
\begin{equation} \label{mag_poten}
A_j \ \in \ L^2_{loc} (\mathbb{R}^3), \ \ j=1,2,3.
\end{equation}
We continue to denote by $S_A$ the Friedrichs extension of
(\ref{MagneticSch}) on $[C_0^{\infty} (\mathbb{R}^3)]^2$.
It is a non-negative self-adjoint operator with no zero modes, and its
form domain $\mathcal{Q}(S_A)$  is the completion of $[C_0^{\infty}
(\mathbb{R}^3)]^2$ with respect to the norm given by
\begin{equation} \label{form_domainSch}
\| \varphi \|_{1,A} \ = \ \biggl\{ \biggl\| \biggl( \frac{1}{i} \bbf{\nabla}
+ \vec{A}\biggr) \varphi \biggr\|^2  + \| \varphi \|^2 \biggr\}^{1/2}.
\end{equation}
The operator realisation of $\mathbb{P}_A$ is given in the first
lemma.
\begin{lemma} \label{oper_realis}  Let \  $|\vec{B}| \in L^{3/2}
(\mathbb{R}^3)$. Then the sesquilinear form
\begin{equation} \label{sesq_form}
p_A [\varphi, \psi] \ = \  (\mathbb{P}_A \varphi, \psi), \ \ \varphi,
\psi \in [C_0^{\infty} (\mathbb{R}^3)]^2
\end{equation}
is symmetric, closable and non-negative in $\mathcal{H}$. The associated
 self-adjoint operator $\mathbb{P}_A$ has form domain $\mathcal{Q}(S_A)$.
\end{lemma}

\begin{proof} Given $\varepsilon > 0$, we may write $|\vec{B}|= B_1 +
B_2$, where $\| B_1 \|_{L^{3/2}(\mathbb{R}^3)} < \varepsilon$  and
 $\| B_2 \|_{L^{\infty}(\mathbb{R}^3)} < C_{\varepsilon}$, for some
constant $C_{\varepsilon}$ depending on $\varepsilon$. Then
\[ (\mathbb{P}_A \ \varphi , \varphi) \ = \ (S_A \ \varphi , \varphi)
+ ((\vec{\sigma} \cdotp \vec{B}) \ \varphi , \varphi)
\]
and
\begin{eqnarray*}
| ((\vec{\sigma} \cdotp \vec{B}) \ \varphi , \varphi) | & \leq& (B_1\
\varphi , \varphi) + (B_2 \ \varphi , \varphi) \\
 &\leq& \| B_1 \|_{L^{3/2}(\mathbb{R}^3)} \| \varphi \|^2_{[L^6
(\mathbb{R}^3)]^2} + C_{\varepsilon} \| \varphi \|^2 \\
&\leq& \varepsilon \gamma^2 \ \| \bbf{\nabla} |\varphi|  \|^2 +  C_{\varepsilon} \| \varphi \|^2
\end{eqnarray*}
by the Sobolev Embedding Theorem, with $\gamma$ the norm of the
embedding $H^1(\mathbb{R}^3)  \hookrightarrow L^6 (\mathbb{R}^3)$,
\[ \ \ \  \leq \varepsilon \gamma^2 \ \biggl\|
\biggl( \frac{1}{i}\bbf{\nabla} + \vec{A}   \biggr)\varphi  \biggr\|^2
+  C_{\varepsilon} \| \varphi \|^2
\]
by the diamagnetic inequality (see \cite[Thm 7.21]{LLB}). The lemma follows from this.
\end{proof}

Hereafter, we shall always assume that
\begin{equation} \label{assum_mag}
|\vec{B} |  \ \in \ L^{3/2} (\mathbb{R}^3).
\end{equation}
The operator
\begin{equation} \label{P}
\mathbb{P} \  := \ \mathbb{P}_A + |\vec{B}|
\end{equation}
may be defined as in Lemma~\ref{oper_realis}, namely, the self-adjoint
operator associated with the form
\begin{equation} \label{form_P}
p[\varphi] \ \equiv \ p[\varphi, \varphi] \ = \ (\mathbb{P}\varphi,
\varphi),
\end{equation}
with form domain $\mathcal{Q} (S_A)$.  As for $S_A$, $\mathbb{P}$ has
no zero modes. Thus  $S_A$ and $\mathbb{P}$ are injective and have
dense domains and ranges in $\mathcal{H}$. Furthermore, $\mathcal{D}
(\mathbb{P}^{1/2})=\mathcal{D} (S_A^{1/2}) =\mathcal{Q}(S_A)$.

    The operator of prime interest is $\mathbb{P}_A$.   We shall
write it as $\mathbb{P}_A = \mathbb{P} - |\vec{B}|$, and then,
initially, proceed along lines which are reminiscent of those
described in \cite{BS} for proving the Cwikel-Lieb-Rosenbljum
inequality for the Schr\"{o}dinger operator.  The problem is
essentially reduced to one for an associated operator of
Birman-Schwinger type. The following spaces feature prominently in the
analysis.
\begin{itemize}
\item \  $H^1_A$ is the completion of $\mathcal{D}(S_A^{1/2})$ with
respect to the norm
\begin{equation} \label{magn_norm}
\| \varphi \|_{H^1_A} \ := \ \| S_A^{1/2} \varphi  \|;
\end{equation}
$H^1_0$ \ has norm $\| \varphi \|_{H^1_0} \ := \ \| \bbf{\nabla} \varphi  \|$.
\item \ $\mathbb{H}^1_B$ is the completion of $\mathcal{D}
(\mathbb{P}^{1/2})$ with respect to the norm
\begin{equation} \label{B_norm}
\| \varphi \|_{\mathbb{H}^1_B} \ := \ \| \mathbb{P}^{1/2} \varphi  \| .
\end{equation}
\end{itemize}

{\sc \underline{Remarks}}
\begin{enumerate}
\item \ $[C_0^{\infty}(\mathbb{R}^3)]^2$ is dense in  $H^1_A$ and
$\mathbb{H}^1_B$.
\item \ The space  $H^1_0$ is not a subspace $\mathcal{H}$. However,
for $\varphi \in [C_0^{\infty}(\mathbb{R}^3)]^2$, the Hardy inequality
\[ \int\limits_{\mathbb{R}^3} \frac{|\varphi (\bbf {x})|^2}{|\bbf
{x}|^2} d \bbf {x} \ \leq \ 4 \int\limits_{\mathbb{R}^3}
|\bbf{\nabla}    \varphi (\bbf {x})   |^2   d \bbf {x}
\]
is valid, and this implies that   $H^1_0$  may be identified with the
function space
\begin{equation} \label{func_space}
H^1_0 \ = \ \bigl\{  u \in [H^1_{loc} (\mathbb{R}^3)] ^2  \ : \
\|  u \|^2_{H^1_0}  + \| \ u / | \cdotp | \ \|^2  < \infty \bigr\}
\end{equation}
and $\| \cdotp \|_{H_0^1}$ is equivalent to the norm defined by
\[
\bigl( \|  u \|^2_{H^1_0}  \  + \ \| \ u / | \cdotp | \ \|^2
\bigr)^{1/2} .
\]
\item \ For the spaces  $H^1_A$ and
$\mathbb{H}^1_B$, which also do not lie in $\mathcal{H}$, we have the
natural embedding
\begin{equation} \label{nat_emb}
\mathbb{H}^1_B \hookrightarrow H^1_A.
\end{equation}
Also, by the diamagnetic inequality, $\varphi \mapsto |\varphi|$ maps
$H^1_A$ continuously into $H^1_0$, which, in turn, is continuously
embedded in $[L^6 (\mathbb{R}^3) ]^2$ by the Sobolev Embedding Theorem.
In fact the spaces in (\ref{nat_emb}) are isomorphic when
(\ref{assum_mag})  is satisfied.
\end{enumerate}

    For a magnetic potential $\vec{A}$ satisfying $|\vec{A}| \in
L^3 (\mathbb{R}^3)$, $H^1_A$ can be shown to be continuously embedded
in $H_0^1$. Such a choice of  $\vec{A}$  is possible in view of the
next lemma which is similar to Theorem~A1 in Appendix~A of \cite{FLL}.

\begin{lemma} \label{lemma2.2}
Let  \ $|\vec{B}| \in L^{3/2}(\mathbb{R}^3)$ and define
\begin{equation} \label{gauss}
\vec{A} (\bbf{x}) \ = \ \frac{1}{4 \pi} \int\limits_{\mathbb{R}^3}
\frac{ (    \bbf{x}-\bbf{y} ) }{|  \bbf{x}-\bbf{y} |^{3}}
 \times \vec{B}
(\bbf{y}) d \bbf{y} .
\end{equation}
Then $|\vec{A}| \in L^3 (\mathbb{R}^3)$, $\mathop{curl} \vec{A} =
\vec{B}$, \  $\mathop{div} \vec{A}=0$ in $\mathcal{D}'$
and
\[ \| \vec{A}  \|_{L^3 (\mathbb{R}^3)} \leq C \ \| \vec{B}
\|_{L^{3/2}  (\mathbb{R}^3)}
\]
for some constant $C$.
\end{lemma}

\begin{proof} The proof is similar to that in \cite{FLL}.
The following formal argument for deriving (\ref{gauss}) is
instructive, and will be helpful for obtaining the analogous result in
$\mathbb{R}^n$ for $n>3$.

The set of Hamiltonian quaternions $\mathbb{H}$ is the unitary
$\mathbb{R}-$algebra generated by the symbols $i,j,k$
with the relations
\[
i^2  =   j^2   =  k^2  =  -1
\]
\[
i j = -j i =k, \ \ j k=-k j =i. \ \ k i = - i k = j.
\]
Multiplication is associative but obviously not commutative.

If we identify a magnetic field $\vec{B}= (B_1 , B_2, B_3)$ and a
magnetic potential  $\vec{A}= (A_1 , A_2, A_3)$ with purely imaginary
quaternionic fields on $\mathbb{R}^3$
\[
\bbf{b} = B_1 (\bbf{x}) i  +  B_2   (\bbf{x}) j +   B_3 (\bbf{x}) k ,
\]
\[
\bbf{a} = A_1 (\bbf{x}) i  +  A_2   (\bbf{x}) j +   A_3 (\bbf{x}) k,
\]
 then the equation
\[
\tilde{D} (\bbf{a}) \ = \ \bbf{b},
\]
where $\tilde{D} = i \frac{\partial}{\partial x_1} +  j
\frac{\partial}{\partial x_2} +  k  \frac{\partial}{\partial x_3}$, is
equivalent to
\[
 \mathop{curl}   \vec{A} = \vec{B}, \ \ \ \mbox{and} \ \ \mathop{div}
 \vec{A} = 0.
\]
We can solve the equation $\tilde{D} (\bbf{a}) \ = \ \bbf{b}$ by the
convolution of $\bbf{b}$ with the Green's function of $\tilde{D}$. Since
$\tilde{D}^2 = - \Delta $ then $\tilde{D} (G(\bbf{x}))$ is the Green's
function for $\tilde{D}$ if $G(\bbf{x})$ is the Green's function
$\frac{1}{4 \pi} \frac{1}{|\bbf{x}|}$ for $- \Delta$.  The identity
(\ref{gauss}) is exactly this convolution of $\bbf{b}$ with $\tilde{D} (G(\bbf{x}))$.
\end{proof}

\begin{lemma} \label{lemma2.3}
Let \  $|\vec{A}| \in L^3 (\mathbb{R}^3)$,  $|\vec{B}| \in
L^{3/2}(\mathbb{R}^3)$. Then \\
(i) for all $f \in \mathbb{H}^1_B$,
\begin{eqnarray}
|  \  (  \   [ \vec{D} \cdotp \vec{A} +  \vec{A} \cdotp \vec{D}  ]f , f )
\ |
& \leq & \  2 \gamma \ \| \vec{A}  \|_{L^3 (\mathbb{R}^3)}
\  \| f  \|^2_{H_0^1}  \label{1} \\
 & \leq &  \ 2 \gamma \ \| \vec{A}  \|_{L^3 (\mathbb{R}^3)}
\  \| f  \|^2_{\mathbb{H}_B^1}, \label{2}
\end{eqnarray}
where $\vec{D} = \frac{1}{i} \bbf{\nabla}$ and $\gamma$ is the norm of
the Sobolev embedding
 $H_0^1\hookrightarrow L^6 (\mathbb{R}^3)$; \\
(ii) \begin{equation} \label{seq_emb}
\mathbb{H}^1_B \hookrightarrow H^1_A \hookrightarrow H^1_0
\hookrightarrow [L^6 (\mathbb{R}^3)]^2 . \end{equation}
\end{lemma}

\begin{proof} \ (i) Let $\varphi \in [C_0^{\infty}
(\mathbb{R}^3)]^2$. Then
\begin{eqnarray*}
\biggl|  \  (  \   [ \vec{D} \cdotp \vec{A} +  \vec{A} \cdotp \vec{D}  ]
\varphi  ,  \varphi )
\ \biggl|
& \ =  \ & \biggr| \    2 \ \mathrm{Re} \sum\limits_{j=1}^3 (A_j
\varphi  , D_j \varphi)    \
\biggl| \\
  & \ \leq  \ &  2 \| \vec{A} \|_{L^3
(\mathbb{R}^3)} \ \|    \varphi  \|_{[L^6 (   \mathbb{R}^3  )]^2}
 \  \| \bbf{\nabla} \varphi      \| \\
& \ \leq  \ & 2 \gamma \ \| \vec{A} \|_{L^3
(\mathbb{R}^3)}   \ \ \| \bbf{\nabla} \varphi      \|^2
\end{eqnarray*}
Thus (\ref{1}) follows by continuity, and this implies  (\ref{2}) once
 (\ref{seq_emb}) is established. \\
(ii) \ Let  $\varphi \in [C_0^{\infty}
(\mathbb{R}^3)]^2$  and $k>1$. Then
\begin{eqnarray*}
k  (S_{\vec{A}} \ \varphi , \varphi) \ &  = &\ (k-1)(- \Delta \varphi ,
\varphi) +
\bigl( \{  - \Delta + k    [\vec{D} \cdotp \vec{A}   +   \vec{A} \cdotp
\vec{D}  ]  + k^2 |\vec{A}|^2    \}  \ \varphi , \varphi   \bigr) \\
 &   & \  + \ \bigl( [k |\vec{A}|^2 -   k^2 |\vec{A}|^2 ] \varphi , \varphi\bigr)
\\
& =  & (k-1) (- \Delta \varphi ,
\varphi) + (S_{k \vec{A}} \ \varphi , \varphi) - (k^2 -k)
\bigl(  |\vec{A}|^2  \varphi , \varphi\bigr) \\
& \geq &  (k-1) (- \Delta \varphi ,
\varphi)  - (k^2 -k)  \bigl(  |\vec{A}|^2  \varphi , \varphi\bigr),
\end{eqnarray*}
whence
\begin{eqnarray*}
(k-1) \|  \bbf{\nabla} \varphi  \|^2 & \leq &
k \| \varphi \|^2_{H^1_A} + k(k-1) \gamma^2 \  \| \vec{A} \|^2_{L^3
(\mathbb{R}^3)} \  \| \  \bbf{\nabla} |\varphi| \  \|^2 \\
& \leq  & \bigl\{  k + k(k-1) \gamma^2 \   \| \vec{A} \|^2_{L^3
(\mathbb{R}^3)}   \bigr\} \ \| \varphi \|^2_{H^1_A}
\end{eqnarray*}
by the diamagnetic inequality. Thus  (\ref{seq_emb}) is established,
and so (\ref{2}).
\end{proof}

\section{A Birman-Schwinger operator}

    Set
\begin{equation} \label{3.1}
p[\varphi] \ : = \ (\mathbb{P}\varphi,
\varphi), \ \ \  b[\varphi] \ : = \ (|\vec{B}| \varphi,
\varphi)
\end{equation}
on $  [C_0^{\infty} (\mathbb{R}^3)]^2$, so that $p_A = p - b$.
From Lemma~\ref{oper_realis} and the remark after  (\ref{form_P}),
 the operators $\mathbb{P}_A$,
$\mathbb{P}$ associated with $p_A$, $p$ respectively have the same
form domain $\mathcal{Q}$, and this is $\mathcal{D} (\mathbb{P}^{1/2})$ with the
graph norm
\begin{equation} \label{graph_norm}
\bigl( \  \|  \mathbb{P}^{1/2} \varphi \|  + \|  \varphi \|^2      \    \bigr)^{1/2}.
\end{equation}
Also  $  [C_0^{\infty} (\mathbb{R}^3)]^2$ is a form core. It follows
that
\begin{equation} \label{3.3}
\mathcal{Q} \ = \ \mathbb{H}^1_B \ \cap \ \mathcal{H}
\end{equation}
with norm (\ref{graph_norm}); the embedding $\mathbb{H}^1_B
\hookrightarrow [L^6 (\mathbb{R}^3)]^2$ guarantees the completeness,
since convergent sequences in $\mathbb{H}^1_B$ therefore converge
pointwise to their limits, almost everywhere.

    From
\[ 0 \leq b[\varphi] \leq p[\varphi]
\]
it follows that there exists a bounded self-adjoint operator
$\mathcal{B}$ on $\mathbb{H}^1_B$ such that
\begin{equation} \label{3.4}
b[\varphi] \ = \ (\mathcal{B} \ \varphi, \varphi )_{\mathbb{H}^1_B},
\ \ \ \varphi \in \mathbb{H}^1_B .
\end{equation}

    For $\varphi \in \mathcal{R}(\mathbb{P}^{1/2})$, the range of
$\mathbb{P}^{1/2}$,
\begin{equation} \label{3.5}
\|   \mathbb{P}^{-1/2} \ \varphi   \|_{\mathbb{H}^1_B} \ = \ \|
\varphi \|
\end{equation}
and hence, since $\mathcal{D} (\mathbb{P}^{1/2}) $ and
$\mathcal{R}(\mathbb{P}^{1/2})$ are dense subspaces of
$\mathbb{H}^1_B$, $\mathcal{H}$ respectively, $\mathbb{P}^{-1/2} $
extends to a unitary map
\begin{equation} \label{3.6}
U \ : \ \mathcal{H} \longrightarrow \mathbb{H}^1_B, \ \ \
U= \mathbb{P}^{-1/2}  \ \ \mbox{on} \ \  \mathcal{R}(\mathbb{P}^{1/2}).
\end{equation}
Define
\begin{equation} \label{3.7}
\mathcal{S} \ := \ |\vec{B}|^{1/2} U \ : \ \mathcal{H} \longrightarrow
\mathcal{H}
\end{equation}
Note that for $u \in \mathbb{H}^1_B$,
\begin{equation} \label{26}
\| \    |\vec{B}|^{1/2} u  \ \|^2   \leq
\| \vec{B}  \|_{L^{3/2} (\mathbb{R}^3)   } \|  u \|^2_{[L^{6} (\mathbb{R}^3)]^2   }\leq
\mathrm{const} \cdotp
\| u  \|^2_{\mathbb{H}^1_B}
\end{equation}
by (\ref{seq_emb}).

\begin{theorem} \label{theorem3.1}
\begin{eqnarray*}
\mathrm{nul } \ \mathbb{P}_A  & =  & \dim \bigl\{  u  : \ \mathcal{B} u =
u, \ u \in   \mathbb{H}^1_B \ \cap \ \mathcal{H}\bigr\} \\
 & \leq  &  \mathrm{nul } \ F,
\end{eqnarray*}
where $F = 1 - \mathcal{S}  \mathcal{S}^{\ast}$.
\end{theorem}

\begin{proof}
Let $u , \varphi \in \mathcal{D} (\mathbb{P}^{1/2})$. Then
\begin{eqnarray*}
p_A[u,\varphi] & = & p[u,\varphi] - b[u,\varphi] \\
& = & (u- \mathcal{B}u , \varphi)_{\mathbb{H}^1_B}.
\end{eqnarray*}
Hence,  $u \in \ker \mathbb{P}_A \subset \mathcal{D}
(\mathbb{P}^{1/2})$ if and only if $\mathcal{B} u = u$ with $u \in
\mathcal{H}$. Moreover, for any $f,g \in \mathcal{H}$
\[
\bigl( \  \mathcal{S} f ,    \mathcal{S} g      \ \bigr) \ = \
\bigl( \  \mathcal{B} \ U  f ,   U  g      \ \bigr)_{\mathbb{H}^1_{B}},
\]
whence
\[
\bigl( \ [\mathcal{S}^{\ast}  \mathcal{S} -1 ] f ,    g      \ \bigr) \ = \
\bigl( \  [\mathcal{B}-1] \ U  f ,   U  g      \ \bigr)_{\mathbb{H}^1_{B}}
\]
The result follows since $\mathrm{nul} \  [  \mathcal{S}^{\ast}
\mathcal{S} -1  ]  =   \mathrm{nul} \   [  \mathcal{S}  \mathcal{S}^{\ast} -1   ] $.
\end{proof}

    The operator $\mathcal{S}  \mathcal{S}^{\ast}$ is of
Birman-Schwinger type.  We have, in terms of  (\ref{3.7})
\begin{equation} \label{3.8}
\mathcal{S}  \mathcal{S}^{\ast} \ = \  |\vec{B}|^{1/2} U^2
|\vec{B}|^{1/2}
\ \ \ \mbox{on} \ \ \mathcal{D} (\mathbb{P}^{1/2})
\end{equation}
and this extends by continuity to a bounded operator on $\mathcal{H}$.
To see (\ref{3.8}), first observe that for $f \in \mathcal{R}
(\mathbb{P}^{1/2})$,  $g \in \mathcal{D}
(\mathbb{P}^{1/2})$
\begin{eqnarray*}
\bigl( \ f , \mathcal{S}^{\ast}  g \ \bigr) & \ = \ & \bigl( \
\mathcal{S} f ,   g \ \bigr) \ = \  \bigl( \
|\vec{B}|^{1/2} U  f ,   g \ \bigr) \\
& \ = \ &  \bigl( \ U f ,  |\vec{B}|^{1/2}    g \ \bigr) \\
& \ = \ &  \bigl( \ \mathbb{P}^{-1/2} f ,  |\vec{B}|^{1/2}    g \ \bigr) ;
\end{eqnarray*}
note that $ |\vec{B}|^{1/2}    g  \in \mathcal{H}$ by (\ref{26})
and since $|\vec{B}| \in L^{3/2} (\mathbb{R}^3)$. Hence $|\vec{B}|^{1/2}    g  \in \mathcal{D}
(\mathbb{P}^{-1/2})$ and $\mathbb{P}^{-1/2} |\vec{B}|^{1/2}    g  =
\mathcal{S}^{\ast}  g$.
In other words
\begin{equation} \label{3.9}
 \mathcal{S}^{\ast} \ = \   \mathbb{P}^{-1/2}  |\vec{B}|^{1/2}
\ \ \ \mbox{on} \ \ \mathcal{D} (\mathbb{P}^{1/2}),
\end{equation}
whence (\ref{3.8}).

\begin{lemma}
$\mathcal{S}  \mathcal{S}^{\ast} $ \ is compact and
\begin{equation} \label{3.10}
\| \mathcal{S} \|^2 \  \leq \ \gamma^2 \ \|      \vec{B}        \|_{L^{3/2}(\mathbb{R}^3)}
\end{equation}
where $\gamma$ is the norm of $H_0^1 \hookrightarrow [L^6 (\mathbb{R}^3)]^2$.
\end{lemma}

\begin{proof}
This is quite standard, but we give the short proof for
completeness. We show that $|\vec{B}|^{1/2} : \mathbb{H}^1_B
\longrightarrow \mathcal{H}$ is compact. Let $\{ \varphi_n \}$ be
a sequence which converges weakly to zero in $\mathbb{H}^1_B$, and
hence in $H_0^1$ by   (\ref{seq_emb}). Then, in particular $\|
\varphi_n \|_{H_0^1} \leq k$, say. Given $\varepsilon > 0$, set
$|\vec{B}| = B_1 + B_2$ where $B_1 \in C_0^{\infty} (\mathbb{R}^3)$
with support $\Omega_{\varepsilon}$ and $B_1 \leq k_{\varepsilon}$
say, and $\|  B_2  \|_{L^{3/2}(\mathbb{R}^3)} < \varepsilon$.
Then
\begin{eqnarray*}
\|  \    |\vec{B}|^{1/2}    \varphi_n   \  \|^2 &\ \leq \ &
 k_{\varepsilon} \| \varphi_n \|^2_{[L^2  (\Omega_{\varepsilon})]^2} +
\ \gamma^2 \ \|   B_2  \|_{L^{3/2}(\mathbb{R}^3)} \ \| \varphi_n
\|^2_{H_0^1} \\
&\ \leq \ & k_{\varepsilon} \| \varphi_n \|^2_{[L^2
(\Omega_{\varepsilon})]^2} + \  \gamma^2 \ \varepsilon \ \| \varphi_n
\|_{\mathbb{H}^1_B} .
\end{eqnarray*}
The first term on the right-hand side tends to zero as $n \rightarrow
\infty$ by the Rellich Theorem. Consequently $|\vec{B}|^{1/2} : \mathbb{H}^1_B
\longrightarrow \mathcal{H}$ is compact  and hence so is
$\mathcal{S}= |\vec{B}|^{1/2} U$.

    The inequality (\ref{3.10}) follows from (\ref{26}).
\end{proof}

\section{The main result}

    For $t \in (0, \infty)$, replace $\vec{A}$ by $t \vec{A}$ and
denote the corresponding operators by $\mathbb{P}_t$, $\mathcal{S}_t$
and $F_t$. It follows from (\ref{3.10}) that
\[
\| \mathcal{S}_t \|^2 \ \leq \ \gamma^2 \ t \  \| \vec{B}
\|_{L^{3/2}(\mathbb{R}^3)}
\longrightarrow 0
\]
as $t \rightarrow 0$. Hence, $F_t = 1- \mathcal{S}_t
\mathcal{S}_t^{\ast}  $  is such that, for some $t_0 > 0$,
\begin{equation} \label{4.1}
\mathrm{nul } \ F_t \ = \ 0, \ \ \ t\in (0, t_0).
\end{equation}
We proceed to prove that $\{  \mathcal{S}_t
\mathcal{S}_t^{\ast}  \}$ is a real analytic family.
\begin{lemma} \label{lemma4.1}
Let $s \in (0, \infty)$ be fixed, and suppose that $| \vec{A}  | \in
L^3 (\mathbb{R}^3)$. Then, there exists a neighbourhood  $N(s)$ of $s$
such that
\begin{equation} \label{4.2}
\frac{1}{t} \mathcal{S}_t \mathcal{S}_t^{\ast} \ = \
\frac{1}{s} \mathcal{S}_s \mathcal{S}_s^{\ast} +
\sum\limits_{n=1}^{\infty} \ (t-s)^n K_n, \ \ \ t \in N(s),
\end{equation}
where the $K_n$ are bounded operators on $\mathcal{H}$.
\end{lemma}

\begin{proof}
For $\psi, \varphi \in \bigl[  C^{\infty}_0 (\mathbb{R}^3)
\bigr]^2$
\[
\bigl( \    \mathbb{P}_t  \psi, \varphi         \ \bigr) \ = \
\bigl( \   \bigl[ \ - \Delta + t (\vec{D} \cdotp \vec{A} + \vec{A} \cdotp \vec{D})
+ t^2 |\vec{A}|^2  + t (\vec{\sigma} \cdotp \vec{B} + |\vec{B}|)
\ \bigr] \psi  ,  \varphi     \ \bigr)
\]
and
\[
\bigl( \  \bigl[ \  \mathbb{P}_t -   \mathbb{P}_s    \ \bigr]
\psi, \varphi         \ \bigr) \ = \
(t-s) \bigl( \ Q \psi, \varphi         \ \bigr)
\]
where
\[
Q  \ = \ \vec{D} \cdotp \vec{A} + \vec{A} \cdotp \vec{D} + (t+s)
|\vec{A}|^2  + \vec{\sigma} \cdotp \vec{B} + |\vec{B}|.
\]
From (\ref{2}) and since
\begin{eqnarray*}
\bigl( \ |\vec{A}|^2 \psi, \psi \ \bigr)
& \ \leq \ &
\| \vec{A}  \|^2_{L^3 (\mathbb{R}^3)} \  \| \psi  \|^2_{[L^6
(\mathbb{R}^3)]^2} \\
& \ \leq \ & \gamma^2 \  \|  \vec{A}  \|^2_{L^3 (\mathbb{R}^3)}
\ \|  \psi  \|^2_{\mathbb{H}^1_B}
\end{eqnarray*}
and
\[
\bigl( \  |\vec{B}| \ \psi , \psi  \ \bigr)  \  \leq \ \gamma^2
\| \vec{B}  \|_{L^{3/2}(\mathbb{R}^3)} \ \| \psi  \|^2_{\mathbb{H}^1_B}
\]
we have
\[
|\bigl(  \  Q  \psi, \psi  \ \bigr)| \ \leq \ c \| \psi  \|^2_{\mathbb{H}^1_B}
\]
for some constant $c$, and so $R= \mathbb{P}^{-1/2}_s Q
\mathbb{P}^{-1/2}_s$ satisfies
\[
|\bigl(  \  R  \psi, \psi  \ \bigr)| \ \leq \ c \| \psi  \|^2
\]
and extends to an operator in $\mathcal{L} (\mathcal{H})$, the space
of bounded linear operators on  $\mathcal{H}$. Thus, there exists a
neighbourhood  $N(s)$ of $s$ such that for $t \in N(s)$
\[
\bigl\{ \ 1+ (t-s) R \ \bigr\}^{-1} \ = \ \sum\limits_{n=0}^{\infty}
 (t-s)^n (-R)^n
\]
in $\mathcal{H}$. It follows that
\[
\mathbb{P}_t \ = \ \mathbb{P}_s + (t-s)Q \ = \
\mathbb{P}_s^{1/2} [1+ (t-s)R] \mathbb{P}_s^{1/2},
\]
 and
\[
|\vec{B}|^{1/2}  \mathbb{P}_t^{-1} |\vec{B}|^{1/2}  \ = \
|\vec{B}|^{1/2}  \mathbb{P}_s^{-1/2} \bigl[ \  1+ (t-s) R     \
\bigr]^{-1}
\mathbb{P}_s^{-1/2}  |\vec{B}|^{1/2} ;
\]
note that $\mathrm{nul} \ \mathbb{P}_t =0$ for any $t$.
For $f \in \mathcal{R} (  \mathbb{P}_s^{1/2}   )$
\begin{eqnarray*}
s \ \|   \ |\vec{B}|^{1/2}  \mathbb{P}_s^{-1/2} f     \  \|^2
& \ = \ &
s \ \bigl( \    |\vec{B}|  \mathbb{P}_s^{-1/2} f ,
\mathbb{P}_s^{-1/2} f       \ \bigr) \\
& \ \leq \  &
\bigl( \      \mathbb{P}_s^{-1/2} f ,
\mathbb{P}_s^{-1/2} f       \ \bigr)_{\mathbb{H}^1_{sB}} \\
& \ = \ & \| f \|^2 .
\end{eqnarray*}
Hence $ |\vec{B}|^{1/2}  \mathbb{P}_s^{-1/2}$, and $\mathbb{P}_s^{-1/2}
|\vec{B}|^{1/2}$, are bounded on $\mathcal{H}$. We may therefore write
\[
|\vec{B}|^{1/2}  \mathbb{P}_t^{-1} |\vec{B}|^{1/2} \ = \
\sum\limits_{n=1}^{\infty} \ (t-s)^n K_n \ \ + \ |\vec{B}|^{1/2}
 \mathbb{P}_s^{-1} |\vec{B}|^{1/2},
\]
where
\[
K_n  \  =  \  |\vec{B}|^{1/2}  \mathbb{P}_s^{-1/2} (- R)^n
\mathbb{P}_s^{-1/2}   |\vec{B}|^{1/2},
\]
and the series lies in $\mathcal{L} (\mathcal{H})$ for $t \in N(s)$.
The preceeding argument implies that with
$T_t = |\vec{B}|^{1/2}  \mathbb{P}_t^{-1/2}$, $|T_t^{\ast}|^2$ has an
extension in $\mathcal{L} (\mathcal{H})$. It follows from (\ref{3.9})
that
$T_t^{\ast} = \mathcal{S}_t^{\ast}$,
and this yields the lemma.
\end{proof}

    We are now in a position to apply the argument of Anghel in
\cite{A}.
For $[a,b] \subset (0, \infty)$, set
\[
d_t \ = \ \mathrm{nul} \ F_t, \ \ \ \ \ \ d_{\min} \ = \ \min\limits_{t \in
[a,b]} \ d_t.
\]

\begin{lemma} \label{lemma4.2}
The map $t \mapsto d_t$  \ is upper semi-continuous.
\end{lemma}

\begin{proof}
The kernel of $F_t$ is finite-dimensional, and we have the orthogonal
decomposition
\[
\mathcal{H } \ = \ \ker F_t \ \oplus \ (\ker F_t)^{\perp} .
\]
With respect to this decomposition, we can represent $F_t$ as
\[
F_t \ = \
\left( \begin{array}{ll}
0 & \ 0 \\
0 & \ D_t
\end{array} \right) ,
\]
where $D_t :  (\ker F_t)^{\perp}  \longrightarrow   (\ker
F_t)^{\perp}$ and
\begin{equation} \label{4.3}
\| D_t \| \ \geq \ c_t \ > \ 0.
\end{equation}
We are required to prove  that,  for any $t$, there exists a
neighbourhood $N(t)$ such that $d_{t'} \leq d_t$ for all $t' \in
N(t)$.
We can write
\[
F_{t'} \ = \
\left( \begin{array}{ll}
L_{t'} & \  M_{t'} \\
M_{t'}^{\ast}  & \ D_{t} + C_{t'}
\end{array} \right)
\]
where $L_{t'}  : \ker F_t \longrightarrow   \ker F_t$, \ \
$ C_{t'} : (\ker F_t)^{\perp}  \longrightarrow   (\ker
F_t)^{\perp}$ are bounded self-adjoint operators and
$ M_{t'} : (\ker F_t)^{\perp}  \longrightarrow \ker F_t$ is bounded.
As $t' \rightarrow t$, we know from Lemma~\ref{lemma4.1} that
$L_{t'}$, $M_{t'}$ and $C_{t'}$ $\rightarrow 0$ in norm. Choose a
neighbourhood $N(t)$ of $t$ such that $\| C_{t'} \| < c_t$ for
$t' \in N(t)$, where $c_t$ is the constant in (\ref{4.3}). Then
$D_{t} + C_{t'}$ is invertible for all $t' \in N(t)$. The operator
\[
\mathcal{A} \ = \
\left( \begin{array}{lr}
\mathbb{I} & \ \  -M_{t'}  (D_{t} + C_{t'})^{-1} \\
0   & \ \ \ (D_{t} + C_{t'})^{-1}
\end{array} \right) ,
\]
where $\mathbb{I}$ is the identity, is a bounded injection on
$\mathcal{H}$, and we have
\begin{equation} \label{4.4}
\mathcal{A} \cdotp  F_{t'} \ = \
\left( \begin{array}{lr}
 L_{t'}   - M_{t'}  (D_{t} + C_{t'})^{-1}   M_{t'}^{\ast}   & \ \ \  0  \\
 (D_{t} + C_{t'})^{-1}   M_{t'}^{\ast}  & \ \  \ \mathbb{I}
\end{array} \right) .
\end{equation}
It follows that
\begin{eqnarray}
d_{t'} & \ = \  &  \mathrm{nul} (\mathcal{A} \cdotp  F_{t'}) \nonumber
\\
& \ = \  & \dim \bigl\{ \
\ker [ L_{t'}   - M_{t'}  (D_{t} + C_{t'})^{-1}   M_{t'}^{\ast} ]
\nonumber  \\
& & \ \ \ \ \ \ \ \cap \ \ker [(D_{t} + C_{t'})^{-1}   M_{t'}^{\ast} ]  \ \bigr\}  \nonumber \\
& \ = \  &  \mathrm{nul} \  [ L_{t'}^2 +  M_{t'}  M_{t'}^{\ast} ]
\label{4.5}  \\
& \ \leq \  & d_t , \nonumber
\end{eqnarray}
whence the lemma.
\end{proof}

\begin{theorem} \label{theorem4.3} \  For any $c \in (0, \infty)$,
$d_t \ = \ 0$,  and hence $\mathrm{nul} \ \mathbb{P}_{tA} = 0$,
 on $[0,c]$ \ except at a finite number of points.
\end{theorem}

\begin{proof}
    We already know from Theorem~\ref{theorem3.1} and
(\ref{4.1}) that $\mathrm{nul} \ \mathbb{P}_{tA}  \leq d_t \ = \ 0$ in
$(0, t_0)$. It is therefore sufficient to prove the theorem for
$[a,c]$, where $0 < a < t_0$.  Define
\begin{eqnarray*}
J  & \ = \ & \bigl\{ \
t \in [a,c] \ : \ \mbox{there exists a neighbourhood} \ \ N_t  \\
& & \ \ \ \mbox{of} \ \ t, \ \ \mbox{such that} \ \  d_{t'} \ = \ 0 \ \
\mbox{in} \ \ N^{'}_t  = N_t  \setminus \{ t \}
\ \bigr\}.
\end{eqnarray*}
The theorem will follow if we prove that $J=[a,c]$, in view of the
compactness of $[a,c]$. We shall prove that $J$ is both open and
closed. Since $a < t_0$, we know that $J \neq \emptyset$.

    It is clear from Lemma~\ref{lemma4.2} that $J$ is open. To
prove that it is closed, let $\{ t_k \}$ be a sequence in $J$ and
$\lim t_k =t$; we may assume that $d_{t_k} = 0$.
In the notation of the proof of Lemma~\ref{lemma4.2}, set
\[
Q_{t'} \ = \  L_{t'}^2 +  M_{t'}  M_{t'}^{\ast} .
\]
Then, from (\ref{4.5})
\begin{eqnarray}
d_t  & \ = \ & \mathrm{rank} \ Q_{t'} + \mathrm{nul} \ Q_{t'} \\
& \ = \ &  \mathrm{rank} \ Q_{t'} + d_{t'} \label{4.6}
\end{eqnarray}
and
\begin{equation} \label{4.7}
d_t \ = \ \mathrm{rank} \ Q_{t_k} .
\end{equation}
If we can prove that $\mathrm{rank} \ Q_{t'} = d_t $ for all
$t'$ in some deleted neighbourhood $N'$ of $t$,  it will follow from
(\ref{4.6}) that $t \in J$, as required.

    Since $\mathrm{rank} \ Q_{t_k} =  d_t $, then any
minor $\mathrm{Min}_{t'}$ of $Q_{t'}$ of order greater than $d_t $ must vanish when
$t' = t_k$.  Hence, since $t' \mapsto \mathrm{Min}_{t'} $ is analytic, $\mathrm{Min}_{t'}
= 0$ in some  neighbourhood $N$ of $t$, and so $\mathrm{rank} \ Q_{t'}
\leq d_t $ in $N$. By (\ref{4.7}) there exist a minor of
$Q_{t'}$ of order $d_t $ which does not vanish on some
subsequence of $\{ t_k \}$, and hence can have a zero only at $t' =t$
within some neighbourhood  $N$ of $t$. Consequently,
$d_t  \geq \mathrm{rank} \ Q_{t'} \geq d_t $ for
$t' \in N'=N \setminus \{ t \}$, and, the theorem is proved.
\end{proof}

\begin{theorem} \label{theorem4.4} \ The set
$\bigl\{ \    \vec{B}:  \  \mathrm{nul} \  \mathbb{P}_{ A} = 0,
\ \mathrm{curl} \ \vec{A} = \vec{B} \ \ \mbox{and} \ \ | \vec{B} | \in L^{3/2}(\mathbb{R}^3)
\ \bigr\} $ contains an open dense subset of
$[L^{3/2}(\mathbb{R}^3)]^3$.
\end{theorem}

\begin{proof}
Let $\mathcal{S}$ in (\ref{3.7}) be denoted by  $\mathcal{S}_B$ and
set $F_B = 1 - \mathcal{S}_B \mathcal{S}_B^{\ast}$. We shall prove
that
\begin{equation} \label{40}
\bigl\{ \    \vec{B}:  \  \mathrm{nul} \   F_B= 0
 \ \ \mbox{and} \ \ | \vec{B} | \in L^{3/2}(\mathbb{R}^3)
\ \bigr\}
\end{equation}
is an open  subset of
$[L^{3/2}(\mathbb{R}^3)]^3$; the theorem will then follow from
Theorem~3.1
since the density of (\ref{40}) is a consequence of
Theorem~\ref{theorem4.3}.

    For $\varepsilon > 0$, let $\vec{B},$  $\vec{B}_0$ be magnetic
fields which satisfy $\|  \vec{B}- \vec{B}_0     \|_{L^{3/2}
(\mathbb{R}^3)} < \varepsilon $. Then, if $\vec{A},$  $\vec{A}_0$  are
the associated vector potentials given in Lemma~\ref{lemma2.2},
$\|  \vec{A} - \vec{A}_0     \|_{L^{3}
(\mathbb{R}^3)} < c \ \varepsilon $ for some $c > 0$. It follows as in
the proof of Lemma~\ref{lemma4.1} that, with $\mathbb{P}= \mathbb{P}_A +
|\vec{B}|$ and $\mathbb{P}_0 = \mathbb{P}_{A_0}  +
|\vec{B}_0|$,
\[
\mathbb{P} - \mathbb{P}_0 \ = \ \mathbb{V},
\]
where, for $\varphi \in [C_0^{\infty} (\mathbb{R}^3) ]^2$,
\begin{eqnarray*}
\bigl( \   \mathbb{V}   \varphi , \varphi   \  \bigr)
& \ = \ & \bigl( \  \bigl[  \vec{D} \cdotp (\vec{A} - \vec{A}_0  ) +
 (\vec{A} - \vec{A}_0  ) \cdotp \vec{D}  \bigr]   \varphi , \varphi  \
\bigr) \\
&& \ + \bigl( \ \bigl[  \   |\vec{A}|^2  - |\vec{A}_0|^2 \
\bigr]    \varphi , \varphi    \ \bigr) \\
&& \ + \bigl( \ \bigl[  \
\vec{\sigma} \cdotp (\vec{B} - \vec{B}_0)
+ | \vec{B} | - | \vec{B}_0 |
\ \bigr]    \varphi , \varphi    \ \bigr)
\end{eqnarray*}
and
\begin{eqnarray*}
\bigl|\  \bigl( \   \mathbb{V}   \varphi , \varphi   \  \bigr) \ \bigr|
& \ \leq \ &  c \  \bigl[ \|   \vec{A} - \vec{A}_0    \|_{L^3 (\mathbb{R}^3)}
+  \|  \vec{B} - \vec{B}_0  \|_{L^{3/2}(\mathbb{R}^3)} \bigr]
\
\| \varphi  \|^2_{H_0^1} \\
& \ \leq \ &  c' \ \varepsilon \ \| \varphi  \|^2_{\mathbb{H}^1_{B_0}}
\end{eqnarray*}
on using Lemmas~\ref{lemma2.2} and \ref{lemma2.3} and H\"{o}lder's
inequality.
Moreover,
$\mathbb{U} = \mathbb{P}_0^{-1/2} \ \mathbb{V} \ \mathbb{P}_0^{-1/2}$
satisfies
\[
\bigl|\  \bigl( \   \mathbb{U}   \varphi , \varphi   \  \bigr) \
\bigr| \
\leq  \  c' \ \varepsilon \ \| \varphi  \|^2
\]
and
\[
|\vec{B}|^{1/2} \ \mathbb{P}^{-1} \ |\vec{B}|^{1/2} \ = \
|\vec{B}|^{1/2} \ \mathbb{P}_0^{-1/2} \ \bigl[   \  1 + \mathbb{U} \
\bigr]^{-1} \ \mathbb{P}_0^{-1/2} \ |\vec{B}|^{1/2}
\]
for $\varepsilon$ sufficiently small. Also,  as $\varepsilon
\rightarrow 0$,
\[
|\vec{B}|^{1/2} \ \mathbb{P}_0^{-1/2} \  -  \ |\vec{B}_0|^{1/2} \
\mathbb{P}_0^{-1/2} \
\longrightarrow  \ 0
\]
in $\mathcal{L} (\mathcal{H})$.
It follows that,  as $\varepsilon
\rightarrow 0$,  \ $F_B \rightarrow  F_{B_0}$ in $\mathcal{L}
(\mathcal{H})$, and that, as in the proof of Lemma~\ref{lemma4.2}, the map
\[
\vec{B} \ \  \longmapsto \ \ \mathrm{nul} \ F_B
\]
is upper semi-continuous. The set (\ref{40}) is therefore open and the
theorem is proved.
\end{proof}




\end{article}
\end{document}